# MEAN VALUES OF DEDEKIND SUMS

J. B. Conrey  
Eric Fransen  
Robert Klein  
Clayton Scott

For a positive integer $k$ and an arbitrary integer $h$, the Dedekind sum $s(h,k)$ is defined by

$$s(h,k) = \sum_{a=1}^{k} \left(\left(\frac{a}{k}\right)\right) \left(\left(\frac{ah}{k}\right)\right)$$

where

$$((x)) = \begin{cases} x - [x] - \frac{1}{2} & \text{if } x \notin Z \\ 0 & \text{if } x \in Z \end{cases}$$

This sum was first studied by Dedekind because of the prominent role it plays in the transformation theory of the Dedekind eta-function, which is a modular form of weight $1/2$ for the full modular group $SL_2(Z)$. There is an extensive literature about the Dedekind sums. Rademacher [8] has written an introductory book on the subject.

The most famous property of the Dedekind sums is the reciprocity formula

$$s(h,k) + s(k,h) = \frac{h^2 + k^2 + 1}{12hk} - \frac{1}{4}$$

for positive coprime integers $h$ and $k$. A three term version of this formula was discovered by Rademacher (see Lemma 1).

Some of the research on the Dedekind sums has involved the question of the distribution in the plane of the points $(h/k, s(h,k))$. (Note: it can be shown that $s(h,k) = s(gh, gk)$ so that $s$ is really a function on the rationals.) It is known that these points are dense in the plane. Recently Vardi [9] has shown that the numbers $s(h,k)/\log k$ with $1 \leq h \leq k$, $1 \leq k \leq K$, and $(h,k) = 1$ have a limiting distribution (the Cauchy distribution) as $K \to \infty$. In particular, for almost all fractions on $[0,1]$ with denominator $\leq K$, it is the case that the Dedekind sum is bounded by $\log^{1+\alpha} K$ for any $\alpha > 0$. In this paper, we are concerned with large values of Dedekind sums. We measure these by considering the $2m$th moments of Dedekind sums averaged over reduced fractions in $[0,1]$ with denominator $k$. For

---

Research of the first author supported in part by a grant from the NSF. Research of the second and third authors supported in part by an REU grant from NSF. Research of the fourth author supported in part by Oklahoma State University.

simplicity, we at first restrict our attention to prime $k$. Later we consider the case of composite $k$.

One should observe from the reciprocity formula that

$$s(\pm 1, k) = \pm \frac{k}{12} + O(1)$$

so that the $2m$th moment of the $s(h, k)$ is trivially at least $2(k/12)^{2m}$. In fact, this is the correct order of magnitude. If one studies a graph of $s(h, k)$ for a large prime $k$, one sees a large positive spike near the origin corresponding to the large value of $s(1, k)$ and a large negative spike near 1 corresponding to the large negative spike at $(k-1)/k$. Further inspection reveals a large positive spike to the right of $1/2$, and a large negative spike slightly to the left of $1/2$. There are other spikes near rational points with small denominators. These observations suggest that it might be possible to investigate the mean values in question via the circle method of Hardy and Littlewood. It turns out that it is possible to give a satisfactory proof of an asymptotic formula for the sum in question exactly along these lines. The proof is quite short, possibly the shortest example of a proof that uses these techniques.

We note that these means have arisen in the literature in connection with mean values of Dirichlet $L$-functions at the point 1: Walum [10] has shown that for prime $k$,

$$\sum_{\substack{\chi \bmod k \\ \chi(-1)=-1}} |L(1,\chi)|^4 = \frac{\pi^4(k-1)}{k^2} \sum_{h=1}^{k} s(h,k)^2.$$

Walum was unable to find an asymptotic for the mean square of the Dedekind sum. However, Heath-Brown [4] notes in his review of Walum's paper that by Dirichlet series techniques one can show that the left hand side of the above is

$$\sim \frac{k}{2} \sum_{n=1}^{\infty} \frac{d(n)^2}{n^2}.$$

Here $d(n)$ is the number of positive divisors of $n$. The sum here is well known to be $\zeta(2)^4/\zeta(4)$.

**Theorem 1.** *Suppose that $k$ is a large prime number. Then,*

$$\sum_{h=1}^{k-1} s(h,k)^{2m} = 2\frac{\zeta(2m)^2}{\zeta(4m)} \left(\frac{k}{12}\right)^{2m} + O\left(\left(k^{9/5} + k^{2m-1+\frac{1}{m+1}}\right) \log^3 k\right).$$

Remark. In an effort to keep the exposition relatively straightforward we have not tried to obtain the best possible error terms here. With more work we can obtain the error term $k^{2m-1+1/(2m)} \log^3 k$ for all $m \geq 1$.

**Lemma 1.** *Suppose that $x$, $y$, and $z$ are positive integers which are pairwise relatively prime. Then*

$$s(x\overline{y}, z) + s(y\overline{z}, x) + s(z\overline{x}, y) = \frac{x^2 + y^2 + z^2}{12xyz} - \frac{1}{4}$$

where the bar over a variable denotes the inverse of that variable modulo the product of the other two variables.

This is Rademacher's famous formula for the Dedekind sums. We remark that the formula holds with one of the variables being negative, provided that we adopt the convention that if $k > 0$, then $s(h, -k) = -s(h, k) - 1/2$.

**Lemma 2.** *Suppose that $k \geq 1$ and that $1 \leq h \leq k$ with $(h, k) = 1$. Suppose that $1 \leq q \leq k$ with $(q, k) = 1$, and that $(a, q) = 1$. Define $\epsilon = hq - ak$. Then*

$$s(h, k) = \frac{k}{12q\epsilon} + O(|s(a, q)| + |\epsilon|).$$

*For $\epsilon = 0$ the main term should be interpreted as 0.*

*Proof.* We apply Lemma 1 with $z = k$, $y = q$, and $x = \epsilon$. Note that

$$x\overline{y} = \epsilon\overline{q} \equiv h \mod k$$

since $hq \equiv \epsilon \mod k$. Also, $k\overline{\epsilon} \equiv -\overline{a} \mod q$. Since $s(-\overline{a}, q) = -s(a, q)$, the result follows.

**Lemma 3.** *(Dirichlet's Theorem) Given any real number $\alpha$ and any $Q > 1$, there exists a positive integer $q \leq Q$ and an integer $a$ such that*

$$\left| \alpha - \frac{a}{q} \right| < \frac{1}{qQ}.$$

See Theorem 36 of Hardy and Wright [3].

**Lemma 4.** *(Barkan [1], Hickerson [5], Knuth [6]) Suppose that $q$ is a positive integer and $1 \leq a \leq q$ with $(a, q) = 1$. Let the continued fraction expansion of $a/q$ be given by*

$$\frac{a}{q} = \langle 0; d_1, d_2, \ldots, d_\ell \rangle.$$

*Finally, define $d$ by $ad \equiv -1 \mod q$ and $1 \leq d \leq q$. Then,*

$$s(a, q) = \begin{cases} \frac{1}{12}\left(\frac{a+d}{q} + d_1 - d_2 \pm \cdots - (-1)^\ell d_\ell\right) - \frac{1}{4} & \text{if } \ell \text{ is odd} \\ \frac{1}{12}\left(\frac{a-d}{q} + d_1 - d_2 \pm \cdots - (-1)^\ell d_\ell\right) & \text{if } \ell \text{ is even} \end{cases}.$$

**Lemma 5.** *(Knuth-Yao[7]) Let $N(a, q)$ be the sum of the partial quotients of the continued fraction expansion of $a/q$. Then*

$$\sum_{\substack{1 \leq a \leq q \\ (a,q)=1}} N(a, q) \ll q \log^2 q.$$

Knuth and Yao actually give an asymptotic formula for this sum, but for our purposes the upper bound suffices.

**Lemma 6.** *As $q \to \infty$ through positive integers,*

$$\sum_{\substack{1 \le a \le q \\ (a,q)=1}} |s(a,q)| \ll q \log^2 q.$$

*Proof.* Using the notation and results of Lemma 4 we have

$$s(a,q) \ll \sum_{j=1}^{\ell} d_j.$$

The result then follows immediately from Lemma 5.

*Proof of Theorem 1.* Let $S$ denote the sum in question. We consider all of the fractions $a/q$ with $a$ and $q$ positive coprime integers with $q \le Q_1$. Here $Q_1 < k$ is a parameter to be chosen later. To each such pair of $a$ and $q$, we associate the open interval $I_{a,q} = (a/q - 1/qQ_1, a/q + 1/qQ_1)$. By Lemma 3, for every $h$ with $1 \le h < k$ the fraction $h/k$ belongs to at least one of the intervals $I_{a,q}$. Moreover, if $Q \le Q_1/2$, then the smaller collection of intervals for which $q \le Q$ have no overlap, i.e. if $q, q' \le Q$, then $I_{a,q} \cap I_{a',q'} = \phi$ for any $a, a'$. To check this assertion, we observe that

$$\frac{1}{qQ_1} + \frac{1}{q'Q_1} = \frac{q+q'}{qq'Q_1} < \frac{2Q}{qq'Q_1} \le \frac{1}{qq'} \le \left|\frac{a}{q} - \frac{a'}{q'}\right|.$$

(We further remark, for later use, that no $a'/q'$ belongs to an $I_{a,q}$ if $a/q \ne a'/q'$. For $a'/q' \in I_{a,q}$ implies that $|a/q - a'/q'| < 1/(qQ_1)$ while $a/q \ne a'/q'$ implies that $|a/q - a'/q'| \ge 1/(qq') \ge 1/(qQ_1)$. A consequence is that no $h/k$ can belong to more than two of the intervals $I_{a,q}$.)

Returning to our argument we now have

$$S = \sum_{q \le Q} \sum_{\substack{a \le q \\ (a,q)=1}} \sum_{\frac{h}{k} \in I_{a,q}} s(h,k)^{2m} + O\left(\sum_{Q < q \le Q_1} \sum_{\substack{a \le q \\ (a,q)=1}} \sum_{\frac{h}{k} \in I_{a,q}} s(h,k)^{2m}\right)$$

Now we apply Lemma 2 to approximate $s(h,k)$. (This is the only place that we make use of the hypothesis that $k$ is prime – to ensure that $(q,k) = 1$.) We use the relation

$$(A+B+C)^{2m} = A^{2m} + O(|A|^{2m-1}(|B|+|C|) + |B|^{2m} + |C|^{2m})$$

which holds for arbitrary $A, B, C$. Then we have $S = M + E_1 + E_2$, where

$$M = \sum_{q \le Q} \sum_{\substack{a \le q \\ (a,q)=1}} \sum_{\frac{h}{k} \in I_{a,q}} \left(\frac{k}{12q\epsilon}\right)^{2m},$$

$$E_1 \ll \sum_{Q<q\leq Q_1} \sum_{\substack{a\leq q \\ (a,q)=1}} \sum_{\frac{h}{k}\in I_{a,q}} \left(\frac{k}{12q\epsilon}\right)^{2m},$$

and

$$E_2 \ll \sum_{q\leq Q_1} \sum_{\substack{a\leq q \\ (a,q)=1}} \sum_{\frac{h}{k}\in I_{a,q}} \left(\left(\frac{k}{q|\epsilon|}\right)^{2m-1}(|s(a,q)|+|\epsilon|) + |s(a,q)|^{2m} + |\epsilon|^{2m}\right).$$

We first estimate $E_2$. Now $2m-1 \geq 1$ and $|\epsilon| \geq 1$. Also, $|s(a,q)| \ll q$ trivially and

$$\epsilon = kq\left(\frac{h}{k}-\frac{a}{q}\right) \ll \frac{k}{Q_1}.$$

Thus,

$$E_2 \ll \sum_{q\leq Q_1} \sum_{\substack{a\leq q \\ (a,q)=1}} \sum_{\frac{h}{k}\in I_{a,q}} \left((k/q)^{2m-1}(|s(a,q)|+1) + q^{2m-1}|s(a,q)| + (k/Q_1)^{2m}\right).$$

The last term here is independent of all of the variables of summation and so it's contribution to the sum is

$$\ll k^{2m+1}/Q_1^{2m}$$

since the trivial estimate for the sum with summand 1 is $2k$ (using the earlier remark that any $h/k$ belongs to at most two intervals $I_{a,q}$).

Now we note that

$$\sum_{\frac{h}{k}\in I_{a,q}} 1 \ll 1 + \frac{k}{qQ_1}.$$

Thus, $E_2 \ll$

$$\sum_{q\leq Q_1} \left(1+\frac{k}{qQ_1}\right) \sum_{\substack{a\leq q \\ (a,q)=1}} \left((k/q)^{2m-1}(|s(a,q)|+1) + q^{2m-1}|s(a,q)|\right) + k^{2m+1}/Q_1^{2m}.$$

By Lemma 6,

$$E_2 \ll \log^2 Q_1 \sum_{q\leq Q_1} \left(1+\frac{k}{qQ_1}\right)\left(\frac{k^{2m-1}}{q^{2m-2}}+q^{2m}\right) + k^{2m+1}/Q_1^{2m}$$

$$\ll \left(k^{2m-1}Q_1^\delta + Q_1^{2m+1} + \frac{k^{2m}}{Q_1} + kQ_1^{2m-1}\right)\log^3 Q_1 + k^{2m+1}/Q_1^{2m}$$

where $\delta=1$ if $m=1$ and $\delta=0$ otherwise.

Next we work on the inner sums over $a$ and $h$ for both $M$ and $E_1$. In particular, let

$$T = \sum_{\substack{a\leq q \\ (a,q)=1}} \sum_{\frac{h}{k}\in I_{a,q}} \epsilon^{-2m}$$

where $\epsilon = hq - ak$. We change the sum over $h$ into a sum over $\epsilon$. The conditions become $\epsilon \equiv -ak \mod q$ and $|\epsilon| \leq k/Q_1$. Thus,

$$T = \sum_{\substack{a \leq q \\ (a,q)=1}} \sum_{\substack{\epsilon \equiv -ak \mod q \\ |\epsilon| \leq \frac{k}{Q_1}}} \epsilon^{-2m}.$$

Now we can carry out the summation over $a$ with the effect that the congruence condition on $\epsilon$ becomes the condition that $\epsilon$ and $q$ are coprime. Thus,

$$T = \sum_{\substack{(\epsilon,q)=1 \\ |\epsilon| \leq \frac{k}{Q_1}}} \epsilon^{-2m}$$

from which we see that

$$T = 2 \sum_{\substack{\epsilon=1 \\ (\epsilon,q)=1}}^{\infty} \epsilon^{-2m} + O\left((Q_1/k)^{2m-1}\right).$$

In particular, $T \ll 1$ so that

$$E_1 \ll k^{2m} \sum_{Q < q \leq Q_1} q^{-2m} \ll k^{2m}/Q^{2m-1}.$$

Now

$$\sum_{\substack{\epsilon=1 \\ (\epsilon,q)=1}}^{\infty} \epsilon^{-2m} = \zeta(2m) \prod_{p|q} (1 - 1/p^{2m})$$

where the product is over primes $p$ which divide $q$. Let

$$F(q) = \prod_{p|q} \left(1 - \frac{1}{p^{2m}}\right).$$

Then we have

$$M = 2\zeta(2m) \left(\frac{k}{12}\right)^{2m} \sum_{q \leq Q} \frac{F(q)}{q^{2m}} + O\left(kQ_1^{2m-1}\right).$$

If we extend the sum over $q$ to infinity, we see that

$$M = 2\zeta(2m) \left(\frac{k}{12}\right)^{2m} \sum_{q=1}^{\infty} \frac{F(q)}{q^{2m}} + O\left(k^{2m}/Q^{2m-1} + kQ_1^{2m-1}\right).$$

Now $F$ is a multiplicative function. Therefore, the sum over $q$ can be expressed as an infinite product. We have

$$\sum_{q=1}^{\infty} \frac{F(q)}{q^{2m}} = \prod_p \left(1 + \sum_{j=1}^{\infty}(1 - 1/p^{2m})/p^{2jm}\right)$$
$$= \prod_p \left(1 + \frac{(1 - 1/p^{2m})}{p^{2m} - 1}\right)$$
$$= \prod_p \left(\frac{1 - 1/p^{4m}}{1 - 1/p^{2m}}\right)$$
$$= \frac{\zeta(2m)}{\zeta(4m)}.$$

Thus, we finally have

$$S = 2\frac{\zeta(2m)^2}{\zeta(4m)} \left(\frac{k}{12}\right)^{2m}$$
$$+ O\left(\left(\frac{k^{2m}}{Q_1} + kQ_1^{2m-1} + \frac{k^{2m}}{Q^{2m-1}} + \frac{k^{2m+1}}{Q_1^{2m}} + Q_1^{2m+1} + k^{2m-1}Q_1^{\delta}\right)\log^3 Q_1\right)$$

where the only stipulation is that $Q \leq Q_1/2 < k/2$. Choosing $Q = Q_1/2$ and $Q_1 = k^{3/5}$ if $m = 1$ or $Q_1 = k^{1-1/(m+1)}$ if $m > 1$ we obtain the Theorem.

Now we turn our attention to the case of composite $k$. Here we encounter difficulties when $(q,k) > 1$. In this situation, Rademacher's formula is no longer applicable. Instead we use a formula recently discovered by Hall and Huxley [2] which extends Rademacher's formula.

**Lemma 7.** *Suppose that $a, b, c, d, h, k$ are positive integers satisfying $ad - bc = 1$ and $(h, k) = 1$. Let*

$$\begin{pmatrix} a & b \\ c & d \end{pmatrix} \begin{pmatrix} h \\ k \end{pmatrix} = \begin{pmatrix} x \\ y \end{pmatrix}.$$

*Then*

$$s(a, c) + s(h, k) - s(x, y) = \frac{c^2 + k^2 + y^2}{12cky} - \frac{1}{4}.$$

This lemma is essentially equation (26) of Hall and Huxley's paper. The analogue of Lemma 2 is

**Lemma 8.** *Suppose that $a, q, h, k > 0$ with $(a, q) = 1$, and $(h, k) = 1$. Suppose further that $\epsilon = hq - ak$ satisfies $|\epsilon| \leq k/q$. Then*

$$s(h, k) = \frac{k}{12q\epsilon} + O(|s(a, q)| + |\epsilon| + 1).$$

*If $\epsilon = 0$, then the main term should be interpreted as 0.*

The main difference between this lemma and Lemma 2 is that in Lemma 2 we required that $(q, k) = 1$, whereas here we do not. The condition that $|\epsilon|$ not be too large is probably not really necessary.

*Proof.* We consider two cases according to whether $\epsilon$ is positive or negative. (If $\epsilon = 0$, the statement is trivial, since in that case $s(h, k) = s(a, q)$.) So suppose first of all that $\epsilon < 0$. Let $b$ and $d$ be positive integers such that $ad - bq = 1$ and $0 < d < q$. Then

$$\begin{pmatrix} d & -b \\ -q & a \end{pmatrix} \begin{pmatrix} h \\ k \end{pmatrix} = \begin{pmatrix} f \\ -\epsilon \end{pmatrix}.$$

This equation may be rewritten as

$$\begin{pmatrix} a & b \\ q & d \end{pmatrix} \begin{pmatrix} f \\ -\epsilon \end{pmatrix} = \begin{pmatrix} h \\ k \end{pmatrix}.$$

In order to apply Lemma 7, we need to verify that $f > 0$. To do this we observe that

$$\frac{f}{kd} = \frac{h}{k} - \frac{b}{d} = \frac{h}{k} - \frac{a}{q} + \frac{1}{qd}$$
$$= \frac{\epsilon}{qk} + \frac{1}{qd} = \frac{1}{q}\left(\frac{\epsilon}{k} + \frac{1}{d}\right)$$
$$> 0$$

since $1/d > 1/q$ and $\epsilon/k \geq -1/q$. Applying Lemma 7, we find that

$$s(a, q) + s(f, -\epsilon) - s(h, k) = \frac{k^2 + q^2 + \epsilon^2}{-12kq\epsilon} - \frac{1}{4}$$

or

$$s(h, k) - s(a, q) - s(f, -\epsilon) = \frac{k^2 + q^2 + \epsilon^2}{12kq\epsilon} + \frac{1}{4}.$$

The assertion of the lemma for this case now follows.

Now we consider the case that $\epsilon > 0$. Let $b$ and $d$ be positive integers with $d < q$ and $ad - bq = 1$. Then

$$\begin{pmatrix} b & -d \\ -a & q \end{pmatrix} \begin{pmatrix} k \\ h \end{pmatrix} = \begin{pmatrix} f \\ \epsilon \end{pmatrix}$$

so that

$$\begin{pmatrix} q & d \\ a & b \end{pmatrix} \begin{pmatrix} f \\ \epsilon \end{pmatrix} = \begin{pmatrix} k \\ h \end{pmatrix}.$$

The argument that $f > 0$ is exactly as before. Then, by Lemma 7,

$$s(q, a) + s(f, \epsilon) - s(k, h) = \frac{a^2 + \epsilon^2 + h^2}{12a\epsilon h} - \frac{1}{4}.$$

Using the reciprocity formula twice,

$$s(h, k) = s(a, q) - s(f, \epsilon) + \frac{h^2 + k^2 + 1}{12hk} - \frac{a^2 + q^2 + 1}{12aq} + \frac{a^2 + \epsilon^2 + h^2}{12a\epsilon h} - \frac{1}{4}$$
$$= s(a, q) - s(f, \epsilon) + \frac{q^2 + \epsilon^2 + k^2}{12kq\epsilon} - \frac{1}{4}.$$

The lemma follows in this case, too.

Using this lemma and the techniques of the proof of Theorem 1 we are able to show

**Theorem 2.** *Suppose that $k$ is large. Then,*

$$\sum_{\substack{h=1 \\ (h,k)=1}}^{k} s(h,k)^{2m} = f_m(k)\left(\frac{k}{12}\right)^{2m} + O\left(\left(k^{9/5} + k^{2m-1+\frac{1}{m+1}}\right)\log^3 k\right)$$

*where*

$$\sum_{k=1}^{\infty} \frac{f_m(k)}{k^s} = 2\frac{\zeta(2m)^2}{\zeta(4m)}\frac{\zeta(s+4m-1)}{\zeta(s+2m)^2}\zeta(s).$$

*Proof of Theorem 2.* The only difference we encounter is in evaluating the sum that arises in the main term. In our present situation we have

$$M = \sum_{q \leq Q} \sum_{\substack{a \leq q \\ (a,q)=1}} \sum_{\substack{\frac{h}{k} \in I_{a,q} \\ (h,k)=1}} \left(\frac{k}{12q\epsilon}\right)^{2m}.$$

As in the earlier case we extend the sum over $h$ or $\epsilon$ from $-\infty$ to $\infty$ and the sum over $q$ to $\infty$. Recall that $\epsilon = hq - ak$. Thus we find that

$$M = \left(\frac{k}{12}\right)^{2m} f_m(k) + O(k^{2m}/Q^{2m-1} + kQ_1^{2m-1})$$

where

$$f_m(k) = \sum_{q=1}^{\infty} \frac{1}{q^{2m}} \sum_{\substack{a=1 \\ (a,q)=1}}^{q} \sum_{\substack{h=-\infty \\ (h,k)=1}}^{\infty} (hq-ak)^{-2m}.$$

To obtain a useful formula for $f_m(k)$ we form the Dirichlet series generating function

$$U(s) = \sum_{k=1}^{\infty} \frac{f_m(k)}{k^s}.$$

We proceed to find an expression for $U(s)$. We remove the coprimality conditions by use of the Möbius relation

$$\sum_{d|n} \mu(n) = \begin{cases} 1 \text{ if } n = 1 \\ 0 \text{ if } n \neq 1 \end{cases}$$

After rearranging the sums, we have

$$U(s) = \sum_{d,e=1}^{\infty} \frac{\mu(d)}{d^{s+2m}} \frac{\mu(e)}{e^{4m}} \sum_{k=1}^{\infty} \frac{1}{k^s} \sum_{q=1}^{\infty} \frac{1}{q^{2m}} \sum_{a=1}^{q} \sum_{h=-\infty}^{\infty} (hq-ak)^{-2m}$$

$$= \frac{1}{\zeta(4m)\zeta(s+2m)} \sum_{k=1}^{\infty} \frac{1}{k^s} \sum_{q=1}^{\infty} \frac{1}{q^{2m}} \sum_{a=1}^{q} \sum_{\substack{n=-\infty \\ n \equiv -ak \bmod q}}^{\infty} n^{-2m}.$$

Now suppose that $(k,q) = g$. Then $g \mid n$ and the inner double sum is

$$= \frac{1}{g^{2m}} \sum_{a=1}^{q} \sum_{\substack{n=-\infty \\ n \equiv -a\frac{k}{g} \mod \frac{q}{g}}}^{\infty} n^{-2m}.$$

Now as $a$ runs through a complete residue system mod $q/g$, so $-ak/g$ runs through a complete residue system mod $q/g$, since $(k/g, q/g) = 1$. As $a$ runs from 1 to $q$, $a$ runs over $g$ complete residue systems mod $q/g$. Thus, the above is

$$= 2g^{1-2m} \sum_{n=1}^{\infty} n^{-2m} = 2g^{1-2m} \zeta(2m).$$

Therefore,

$$U(s) = \frac{2\zeta(2m)}{\zeta(4m)\zeta(s+2m)} \sum_{k=1}^{\infty} \sum_{q=1}^{\infty} \frac{g^{1-2m}}{k^s q^m}.$$

Now,

$$g^{1-2m} = \sum_{\substack{d \mid q \\ d \mid k}} \sum_{e \mid d} \mu(e) \left(\frac{d}{e}\right)^{1-2m}.$$

So the inner double sum in the last expression for $U(s)$ above is

$$= \sum_{d=1}^{\infty} \frac{\mu(e) \left(\frac{d}{e}\right)^{1-2m}}{d^{s+2m}} \sum_{k=1}^{\infty} \frac{1}{k^s} \sum_{q=1}^{\infty} \frac{1}{q^{2m}}$$
$$= \frac{\zeta(s)\zeta(2m)\zeta(s+4m-1)}{\zeta(s+2m)}.$$

Thus, we finally have

$$U(s) = \frac{2\zeta(2m)^2 \zeta(s) \zeta(s+4m-1)}{\zeta(4m)\zeta(s+2m)^2}.$$

It might be of interest to determine the poles of

$$\sum_{k=1}^{\infty} k^{-s} \sum_{\substack{1 \le h \le k \\ (h,k)=1}} s(h,k)^{2m}.$$

We would like to thank Professor Martin Huxley for pointing out that Lemma 2 could be proven using Lemma 1, and similarly for Lemmas 7 and 8. Previously we had an adhoc method for proving Lemma 2 based solely on the reciprocity formula. We would also like to thank Ilan Vardi and the referee for helpful remarks.

Department of Mathematics, Oklahoma State University, Stillwater, Oklahoma 74078

Department of Mathematics, Oklahoma State University, Stillwater, Oklahoma 74078

Department of Mathematics, Oklahoma State University, Stillwater, Oklahoma 74078

Oklahoma School of Science and Mathematics, 1515 N. Lincoln Blvd., Oklahoma City, 73104